\documentclass[12pt,a4paper]{article}
\usepackage[english]{babel}
\usepackage{amssymb}
\usepackage{inputenc}
\usepackage{amssymb}
\begin{document}

\author{Shavkat A. Ayupov  and   Karimbergen K. Kudaybergenov }

\title{\bf Additive derivations  on algebras of measurable operators }

\maketitle
\begin{abstract}
Given a von Neumann algebra $M$ we  introduce so called central
extension $mix(M)$ of $M$. We show that
 $mix(M)$ is a *-subalgebra in the algebra $LS(M)$
  of all locally measurable operators with respect to $M,$ and this algebra coincides with
 $LS(M)$ if and only if
 $M $ does not admit type II direct summands. We prove that if  $M$ is a properly infinite  von
Neumann algebra then
every additive derivation on the algebra $mix(M)$ is inner. This implies that on
 the algebra $LS(M)$, where $M$ is a type I$_\infty$ or a type III von Neumann algebra, all additive derivations are inner derivations.
\end{abstract}

\medskip



\newpage

\begin{center}
{\bf Introduction}
\end{center}

The present paper continues the series
 of papers \cite{Alb1}-\cite{Alb2} devoted to the investigation of derivations on the
algebra $LS(M)$ of locally measurable  operators affiliated with a
von Neumann algebra $M$ and on its various subalgebras.

Let  $\mathcal{A}$ be an algebra over the field complex number. A linear (additive) operator
$D:\mathcal{A}\rightarrow \mathcal{A}$ is called \emph{a linear (additive) derivation} if it satisfies the
identity  $D(xy)=D(x)y+xD(y)$ for all  $x, y\in \mathcal{A}$ (Leibniz rule). Further linear derivations are called simply \emph{derivations}.
Each element  $a\in \mathcal{A}$ defines a linear derivation  $D_a$ on $\mathcal{A}$ given as
$D_a(x)=ax-xa,\,x\in \mathcal{A}.$ Such derivations $D_a$ are said to be
\emph{inner derivations}. If the element  $a$ implementing the
derivation   $D_a$ on $\mathcal{A},$ belongs to a larger algebra  $\mathcal{B}$
containing  $\mathcal{A}$ (as a proper ideal as usual) then $D_a$ is called
a  \emph{spatial derivation}.

One of the main problems in the theory of derivations is to prove
the automatic continuity, innerness or spatialness of
derivations  or to show the existence of non inner
and discontinuous derivations on various topological algebras.

On this way  A.~F.~Ber, F.~A.~Sukochev, V.~I.~Chilin~\cite{Ber}
obtained necessary and sufficient conditions for the existence of
non trivial derivations on commutative regular algebras. In
particular they have proved that the algebra  $L^{0}(0, 1)$ of all
(classes of equivalence of) complex measurable functions on  the
interval $(0, 1)$ admits non trivial derivations. Independently
A.~G.~Kusraev~\cite{Kus1} by means of Boolean-valued analysis has
also proved the existence of non trivial derivations and
automorphisms on $L^{0}(0, 1).$ It is clear that these derivations
are discontinuous in the measure topology, and therefore they are
neither inner nor spatial. The present authors have conjectured
that the existence of such pathological examples of derivations
deeply depends on the commutativity of the underlying von Neumann
algebra $M.$ In this connection we have initiated the study of the
above problems in the non commutative case
\cite{Alb1}-\cite{Alb2},
 by considering derivations on the algebra $LS(M)$ of
all locally measurable operators affiliated with a von Neumann
algebra $M$ and on various subalgebras of $LS(M).$ In \cite{Alb3}
we have proved
  that every derivation on so called non commutative Arens
  algebras affiliated with an arbitrary von Neumann algebra and a
  faithful normal semi-finite trace is spatial and if the trace
  is finite then all derivations on this algebra are inner. In  \cite{Alb1}
and  \cite{Alb2}    we have
 proved the above mentioned conjecture concerning derivations on
 $LS(M)$ for type I von Neumann algebras.

  Recently this conjecture was also independently confirmed for the type I
case in the paper of A.F. Ber, B. de Pagter and A.F. Sukochev
\cite{Ber1} by means of a representation of measurable operators
as operator valued functions. Another approach to similar problems
in the framework of type I $AW^{*}$-algebras has been outlined in
the paper of A.F. Gutman, A.G.Kusraev and  S.S. Kutateladze
\cite{Gut}.

In \cite{Alb2}  we considered derivations on the algebra $LS(M)$
of all locally measurable operators affiliated with a type I von
Neumann algebra $M$, and also on its subalgebras $S(M)$ -- of
measurable operators and $S(M, \tau)$ of $\tau$-measurable
operators, where $\tau$ is a faithful normal semi-finite trace on
$M.$ We proved that an arbitrary derivation $D$ on each of these
algebras can be uniquely decomposed into the sum $D=D_a+D_\delta$
where the derivation $D_a$ is inner (for $LS(M),$ $S(M)$ and $S(M,
\tau)$) while the derivation $D_\delta$ is an extension of a derivation $\delta$ on the center of the corresponding
algebra.

In the present paper we consider additive derivations on the algebra $LS(M),$
were $M$ is a properly infinite von Neumann algebras.

In section 1 we  introduce the so called  central extension
$mix(M)$ of a von Neumann algebra $M.$ We show that
 $mix(M)$ is a *-subalgebra in the algebra $LS(M)$ and
 this algebra coincides with whole $LS(M)$ if and only if
 $M$ does not contain a direct summand of type II. The
 center $Z(M)$ of $M$ is an abelian von Neumann algebra and
hence it
 is *-isomorphic to $ L^{\infty}(\Omega, \Sigma, \mu)$
   for an appropriate measure space $(\Omega, \Sigma, \mu)$. Therefore the algebra $LS(Z(M))= S(Z(M))$ can be
 identified with the ring $ L^{0}(\Omega, \Sigma, \mu)$ of all
 measurable functions on $(\Omega, \Sigma, \mu)$.
 We also show that $mix(M)$
is a $C^*$-algebra over the ring $S(Z(M))\cong L^{0}(\Omega,
\Sigma, \mu)$ in the sense of \cite{AK}.

In section 2 we give some necessary properties of the topology of convergence locally in measure on $LS(M).$

In section 3  additive derivations on the algebra $mix(M)$ are investigated. We prove that if  $M$ is a properly infinite  von
Neumann algebra then every additive derivation on the algebra
$mix(M)$ is inner. This implies in particular that every additive derivation on the
algebra $LS(M),$ where $M$ is of type I$_\infty$ or of type III, is in fact an inner derivation.
The latter result generalizes Theorem 2.7 from  \cite{Alb2} to
additive derivations and extends it also for type III von Neumann
algebras.

\begin{center}
\textbf{1. Locally measurable operators affiliated with von Neumann algebras}
\end{center}

Let  $H$ be a complex Hilbert space and let  $B(H)$ be the algebra
of all bounded linear operators on   $H.$ Consider a von Neumann
algebra $M$  in $B(H)$ with the operator norm $\|\cdot\|_M.$ Denote by
$P(M)$ the lattice of projections in $M.$

A linear subspace  $\mathcal{D}$ in  $H$ is said to be
\emph{affiliated} with  $M$ (denoted as  $\mathcal{D}\eta M$), if
$u(\mathcal{D})\subset \mathcal{D}$ for every unitary  $u$ from
the commutant
$$M'=\{y\in B(H):xy=yx, \,\forall x\in M\}$$ of the von Neumann algebra $M.$

A linear operator  $x$ on  $H$ with the domain  $\mathcal{D}(x)$
is said to be \emph{affiliated} with  $M$ (denoted as  $x\eta M$) if
$\mathcal{D}(x)\eta M$ and $u(x(\xi))=x(u(\xi))$
 for all  $\xi\in
\mathcal{D}(x).$

A linear subspace $\mathcal{D}$ in $H$ is said to be \emph{strongly
dense} in  $H$ with respect to the von Neumann algebra  $M,$ if

1) $\mathcal{D}\eta M;$

2) there exists a sequence of projections
$\{p_n\}_{n=1}^{\infty}$ in $P(M)$  such that
$p_n\uparrow\textbf{1},$ $p_n(H)\subset \mathcal{D}$ and
$p^{\perp}_n=\textbf{1}-p_n$ is finite in  $M$ for all
$n\in\mathbb{N},$ where $\textbf{1}$ is the identity in $M.$

A closed linear operator  $x$ acting in the Hilbert space $H$ is said to be
\emph{measurable} with respect to the von Neumann algebra  $M,$ if
 $x\eta M$ and $\mathcal{D}(x)$ is strongly dense in  $H.$ Denote by
 $S(M)$ the set of all measurable operators with respect to
 $M.$

A closed linear operator $x$ in  $H$  is said to
be \emph{locally measurable} with respect to the von Neumann
algebra $M,$ if $x\eta M$ and there exists a sequence
$\{z_n\}_{n=1}^{\infty}$ of central projections in $M$ such that
$z_n\uparrow\textbf{1}$ and $z_nx \in S(M)$ for all
$n\in\mathbb{N}.$

It is well-known \cite{Mur} that the set $LS(M)$ of all locally
measurable operators affiliated with  $M$ is a unital *-algebra
when equipped with the algebraic operations of strong addition and
multiplication and taking the adjoint of an operator, and contains
$S(M)$ as a solid *-subalgebra.

 Let   $\tau$ be a faithful normal semi-finite trace on $M.$ We recall that a closed linear operator
  $x$ is said to be  $\tau$\emph{-measurable} with respect to the von Neumann algebra
   $M,$ if  $x\eta M$ and   $\mathcal{D}(x)$ is
  $\tau$-dense in  $H,$ i.e. $\mathcal{D}(x)\eta M$ and given   $\varepsilon>0$
  there exists a projection   $p\in M$ such that   $p(H)\subset\mathcal{D}(x)$
  and $\tau(p^{\perp})<\varepsilon.$
   The set $S(M,\tau)$ of all   $\tau$-measurable operators with respect to  $M$
    is a solid  *-subalgebra in $S(M)$  (see \cite{Nel}).

    Consider the topology  $t_{\tau}$ of convergence in measure or \emph{measure topology}
    on $S(M, \tau),$ which is defined by
 the following neighborhoods of zero:
$$V(\varepsilon, \delta)=\{x\in S(M, \tau): \exists e\in P(M), \tau(e^{\perp})\leq\delta, xe\in
M,  \|xe\|_{M}\leq\varepsilon\},$$  where $\varepsilon, \delta$
are positive numbers, and $\|.\|_{M}$ denotes the operator norm on
$M$.

 It is well-known
\cite{Nel} that $S(M, \tau)$ equipped with the measure topology is a
complete metrizable topological *-algebra.

Note that if the trace $\tau$ is  finite then
$$S(M, \tau)=S(M)=LS(M).$$

Given any family  $\{z_i\}_{i\in I}$ of mutually orthogonal
central projections in $M$ with $\bigvee\limits_{i\in
I}z_i=\textbf{1}$ and a  family of elements $\{x_i\}_{i\in I}$ in
$LS(M)$ there exists a unique element $x\in LS(M)$ such that $z_i
x=z_i x_i$ for all $i\in I.$ This element is denoted by
$x=\sum\limits_{i\in I}z_i x_i$ and it is called \emph{ the
mixing} of $\{x_i\}_{i\in I}$ with respect to $\{z_i\}_{i\in I}$
(see Proposition 1.1 and further remarks in \cite{Alb2}).

By  $mix(M)$ we denote the set of all elements  $x$ from  $LS(M)$ for which there exists a sequence of
mutually orthogonal central projections  $\{z_i\}_{i\in I}$ in  $M$ with $\bigvee\limits_{i\in I}z_i=\textbf{1},$
such that $z_i x\in M$ for all $i\in I,$ i.e.
 $$mix(M)=\{x\in LS(M): \exists z_i\in P(Z(M)), z_iz_j=0, i\neq j, \bigvee\limits_{i\in I}z_i=\textbf{1},
 z_i x\in M, i\in I\},$$
where $Z(M)$ is the center of $M.$ In other words $mix(M)$ is the set of all mixings obtained by families
$\{x_i\}_{i\in I}$ taken from $M.$

 \textbf{Proposition 1.1.} \emph{Let} $M$ -- \emph{be a von Neumann algebras with the center $Z(M).$
 Then}

i)  $mix(M)$ \emph{is a *-subalgebra in} $LS(M)$ \emph{with the
center} $S(Z(M)),$ \emph{where} $S(Z(M))$ \emph{is the algebra of
operators measurable with respect to}  $Z(M)$;

ii) $LS(M)=mix(M)$ \emph{if and only if} $M$ \emph{does not have
direct summands of type II}.

 Proof. i)  It is clear from the definition that  $mix(M)$ is a *-subalgebra in $LS(M)$ and that
its center  $Z(mix(M))$  is contained in  $S(Z(M))=Z(LS(M)).$

 Let us show the converse inclusion. Take
$x\in S(Z(M))$ and let $|x|=\int\limits_{0}^{\infty}\lambda\,de_{\lambda}$
  be the spectral resolution of $|x|.$
 Set
 $$z_1=e_1\quad\mbox{and}\quad z_n=e_n-e_{n-1},\,n\geq2.$$
 Then it clear that $\{z_n\}_{n\in\mathbb{N}}$   is a sequence of mutually orthogonal central projections in  $M$
 such that
 $\bigvee\limits_{n\geq1}z_n=\textbf{1}$  and $z_n x\in Z(M)$
 for all
 $n\in\mathbb{N}.$
Therefore  $x\in mix(M).$ Since
 $x$  commutes with each element from  $LS(M)\supset mix(M),$
 we have that  $x\in Z(mix(M)).$ Thus
 $Z(mix(M))=S(Z(M)).$

ii)  If  $M$ is of type I, then by  \cite[Proposition 1.6]{Alb2} we have   $LS(M)=mix(M).$

Let $M$   have type  III. Since any nonzero projection in $M$ is
infinite
 it follows that $S(M)=M.$  Hence by the definitions of the algebras  $LS(M)$ and  $mix(M)$
we obtain that  $LS(M)=mix(M).$ Thus if $M=N\oplus K$ where $N$ is
a type I and $K$ is a type III von Neumann algebras, i.e. if $M$
does not have type II direct summands, then $LS(M)=mix(M).$

To prove the converse suppose that $M$ is a type II von Neumann
algebra. First assume that $M$ is a type II$_1$  von Neumann algebra with a
faithful normal tracial state $\tau$. Let  $\Phi$ be the
canonical center-valued trace on $M.$

Since   $M$ is of type  II, then  there exists a projection
$p_1\in M$ such that  $$p_1\sim \textbf{1}-p_1.$$ Then
$\Phi(p_1)=\Phi(p_1^{\perp}).$ From
$\Phi(p_1)+\Phi(p_1^{\perp})=\Phi(\textbf{1})=\textbf{1}$ it follows that
$$\Phi(p_1)=\Phi(p_1^{\perp})=\frac{\textstyle 1}{\textstyle 2}\textbf{1}.$$

Suppose that there exist  mutually orthogonal projections
$p_1,\,p_2,\cdots,p_n$ in $M$ such that
$$\Phi(p_k)=\frac{\textstyle 1}{\textstyle 2^{k}}\textbf{1},\,k=\overline{1, n}.$$
Set  $e_n=\sum\limits_{k=1}^{n}p_k.$ Then
$\Phi(e_n^{\perp})=\frac{\textstyle 1}{\textstyle
2^{n}}\textbf{1}.$ Take a projection  $p_{n+1}<e_n^{\perp}$ such that
$$p_{n+1}\sim e_n^{\perp}-p_{n+1}.$$
Then
$$\Phi(p_{n+1})=\frac{\textstyle 1}{\textstyle
2^{n+1}}.$$

Hence there exists a sequence a mutually orthogonal
projections $\{p_n\}_{n\in\mathbb{N}}$ in $M$ such that
$$\Phi(p_n)=\frac{\textstyle
1}{\textstyle 2^{n}}\textbf{1},\,n\in\mathbb{N}.$$
Note that  $\tau(p_n)=\frac{\textstyle 1}{\textstyle 2^{n}}.$
Indeed
$$\tau(p_n)=\tau(\Phi(p_n))=\tau(\frac{\textstyle
1}{\textstyle 2^{n}}\textbf{1})=\frac{\textstyle 1}{\textstyle
2^{n}}.$$

Since
$$\sum\limits_{n=1}^{\infty}n\tau(p_n)=
\sum\limits_{n=1}^{\infty}\frac{n}{2^{n}}<+\infty$$
it follows that  the series
$$\sum\limits_{n=1}^{\infty}np_n$$ converges in measure in  $S(M, \tau).$
Therefore
there exists  $x=\sum\limits_{n=1}^{\infty}np_n\in S(M, \tau).$

Let us  show that $x\in LS(M)\setminus mix(M).$ Suppose that  $zx\in M,$ where  $z$ is a nonzero central
projection. Since any  $p_n$ is a faithful projection we have that  $z p_n\neq 0$ for all $n.$ Thus
$$||zx||_{M}=1||zx||_{M}1=||p_n||_{M}\cdot||zx||_{M}\cdot||p_n||_{M}\geq||zp_nxp_n||_{M}=||zp_n n||_{M}=n,$$
i.e.
$$||zx||_{M}\geq n$$
for all $n\in\mathbb{N}.$ From this  contradiction it follows that
 $x\in LS(M)\setminus mix(M).$

  For a general type II von Neumann
 algebra $M$ take a non zero finite projection $e\in M$ and
 consider the finite type II von Neumann algebra $eMe$ which
 admits a separating family of normal tracial states. Now if $f\in
 eMe$ is the support projection of some tracial state $\tau$ on
 $eMe$ then $fMf$ is a type II$_1$ von Neumann algebra with a
 faithful normal tracial state.Hence as above one can
 construct an element $x\in LS(M)\setminus mix(M).$
 Therefore if $LS(M)=mix(M)$ then $M$ can not have a direct summand of the type II.
 The proof is complete. $\blacksquare$

\medskip

\medskip

\textbf{Remark.} A similar notion (i.e the algebra $mix(A)$) for
arbitrary *-subalgebras $A\subset LS(M)$  was independently
introduced recently by M.A. Muratov and V.I. Chilin \cite{Mur1}.
They denote this algebra by $E(A)$ and called it the central
extension of $A$. In particular if $A=M$ we have $E(M)=mix(M).$
Therefore following \cite{Mur1} we shall say that $mix(M)$ is
\emph{the central extension of} $M$.

An alternative proof of Proposition 1.1 follows also from
  Proposition 2, Theorem 1 and Theorem 3 in \cite{Mur1}.

Let $(\Omega,\Sigma,\mu)$  be a measure space and from now on
suppose
 that the measure $\mu$ has the  direct sum property, i. e. there is a family
 $\{\Omega_{i}\}_{i\in
J}\subset\Sigma,$ $0<\mu(\Omega_{i})<\infty,\,i\in J,$ such that
for any $A\in\Sigma,\,\mu(A)<\infty,$ there exist a countable
subset $J_{0 }\subset J$ and a set  $B$ with zero measure such
that $A=\bigcup\limits_{i\in J_{0}}(A\cap \Omega_{i})\cup B.$

 We denote by  $L^{0}(\Omega, \Sigma, \mu)$ the algebra of all
(equivalence classes of) complex measurable functions on $(\Omega,
\Sigma, \mu)$ equipped with the topology of convergence in
measure.

Consider the algebra  $S(Z(M))$  of operators measurable with
respect to the  center $Z(M)$ of the von Neumann algebra $M.$
Since  $Z(M)$ is an abelian von Neumann algebra  it
 is *-isomorphic to $ L^{\infty}(\Omega, \Sigma, \mu)$
   for an appropriate measure space $(\Omega, \Sigma, \mu)$. Therefore the algebra  $S(Z(M))$ can be
 identified with the algebra $ L^{0}(\Omega, \Sigma, \mu)$ of all
 measurable functions on $(\Omega, \Sigma, \mu)$.

\textbf{Proposition  1.2.} \emph{For any  $x\in mix(M)$ there exists  $f\in S(Z(M))$ such that  $|x|\leq f.$}

Proof. Let  $x=\sum\limits_{i\in I}z_i x\in mix(M), z_i x\in M$ for all $i\in I.$
Put
$$
f=\sum\limits_{i\in I}z_i ||z_ix||_{M}\in S(Z(M)).
$$ Then
$$|x|=|\sum\limits_{i\in I}z_i x|=\sum\limits_{i\in I}z_i |z_ix|\leq\sum\limits_{i\in I}z_i ||z_ix||_{M}=f.$$
The proof is complete. $\blacksquare$

Proposition 1.2 implies that for any $x\in mix(M)$ there exists the following vector-valued norm
 \begin{equation}
 ||x||=\inf\{f\in S(Z(M)): |x|\leq f\}.
\end{equation}

\medskip

\medskip
\bigskip

By the definition we obtain that:

 a) $|x|\leq ||x||$ for all $x\in mix(M);$

 b) if  $x\in mix(M)$ then
  $$||x||=\inf\{f\in S(Z(M)): f\geq 0, f^{-1}x\in M, ||f^{-1}x||_{M}\leq 1\};$$

c) if  $z\in M$ is a central projection then
  $$||zx||=z||x||;$$

 d) if $x\in M$  then  $$||x||_{M}=||||x||||_{M}.$$

\textbf{Proposition 1.3.} \emph{Let  $x\in M.$  Then}  $||x||=\textbf{1}$ \emph{if and only if
$||zx||_{M}=1$ for each nonzero central projection}  $z\in M.$

Proof.  Let  $x\in M,$ $||x||=\textbf{1}.$ Then
$||zx||=z||x||=z$ for each  nonzero central projection  $z\in M.$ Thus
$$||zx||_{M}=||||zx||||_{M}=||z||_{M}=1.$$

Now let  $||zx||_{M}=1$ for each  nonzero central projection  $z\in M,$ and
in particular,  $||x||_{M}=1.$ Thus  $||x||\leq \textbf{1}.$  Suppose that
 $||x||\neq\textbf{1}.$
Then there exist a nonzero  central projection
 $z\in M$ and a number  $0<\varepsilon<1$ such that
 $z||x||\leq \varepsilon z.$ Thus
 $$ ||zx||_{M}\leq \varepsilon ||z||_{M}=\varepsilon<1, $$
 and this  contradicts
   to the equality  $||zx||_{M}=1.$ Hence  $||x||=\textbf{1}.$
   The proof is complete.  $\blacksquare$

A complex linear space  $E$ is said to be normed  by $L^{0}(\Omega, \Sigma, \mu)$ if
there is
 a map  $\|\cdot\|:E\longrightarrow L^{0}(\Omega, \Sigma, \mu)$ such that for any  $x,y\in E, \lambda\in
\mathbb{C},$ the following conditions are fulfilled:

 1) $\|x\|\geq 0; \|x\|=0\Longleftrightarrow x=0$;

 2) $\|\lambda x\|=|\lambda|\|x\|$;

 3)  $\|x+y\|\leq\|x\|+\|y\|$.

 The  pair  $(E,\|\cdot\|)$ is called a lattice-normed  space over $L^{0}(\Omega, \Sigma, \mu).$
A lattice-normed  space  $E$ is called  $d$-decomposable, if for
any $x\in E$ with  $\|x\|=\lambda_{1}+\lambda_{2},$  $
\lambda_{1}, \lambda_{2} \in L^{0}(\Omega, \Sigma,
\mu),\,\lambda_{1}\lambda_{2}=0, \lambda_1, \lambda_2\geq0,$ there
exist $x_{1}, x_{2}\in E$ such that $x=x_{1}+x_{2}$ and
$\|x_{i}\|=\lambda_{i},\,i=1,2$.

A net $(x_{\alpha})$ in $E$ is said to be $(bo)$-convergent to
$x\in E$, if the net $\{\|x_{\alpha}-x\|\}$ $(o)$-converges (i.e.
almost everywhere converges) to zero in
 $L^{0}(\Omega, \Sigma, \mu).$

 A lattice-normed  space $E$ which is $d$-decomposable and  complete with respect to the  $(bo)$-convergence
  is called a
 \emph{Banach -- Kantorovich space}.

 It is known  that every Banach -- Kantorovich space $E$ over
 $L^{0}(\Omega, \Sigma, \mu)$ is a module over $L^{0}(\Omega, \Sigma, \mu)$ and  $\|\lambda x\|=|\lambda|\|x\|$
 for all $\lambda\in L^{0}(\Omega, \Sigma, \mu),\, x\in E$ (see \cite{Kusr}).

 Let $\mathcal{A}$ be an arbitrary Banach -- Kantorovich space
 over $L^{0}(\Omega, \Sigma, \mu)$ and let $\mathcal{A}$ be a
 *-algebra such that $(\lambda x)^{\ast}=\overline{\lambda} x^{\ast},$ $(\lambda x)y=\lambda (xy)=x(\lambda y)$
for all $\lambda\in L^{0}(\Omega, \Sigma, \mu),$ $x, y\in \mathcal{A}.$ $\mathcal{A}$ is called a $C^\ast$-\emph{algebra
over }$L^{0}(\Omega, \Sigma, \mu)$
 if $||x y||\leq ||x||||y||,$ $||xx^{\ast}||=||x||^2$ for all $x, y \in \mathcal{A}$ (see \cite{AK}).

The main result of this section is the following.

 \textbf{Proposition  1.4.} \emph{Let} $M$ \emph{be a von Neumann algebra with the center}
 $Z\cong L^{\infty}(\Omega, \Sigma, \mu)$ \emph{and let }
  $||\cdot||$ \emph{be the} $S(Z(M))$-\emph{valued
 norm on} $mix(M)$ \emph{defined by} (1).
  \emph{Then }
 $(mix(M), ||\cdot||)$ \emph{is a $C^{\ast}$-algebra over } $S(Z(M))\cong L^{0}(\Omega, \Sigma, \mu).$

Proof. Let  $x\in mix(M), x\neq 0$ and let
   $|x|=\int\limits_{0}^{\infty}\lambda\,de_{\lambda}$ be the spectral resolution  of $|x|.$
Then there exists
 $\lambda_0>0$ such that  $e_{\lambda_0}\neq 0.$ Take an element  $f\in S(Z(M))$ such that
 $|x|\leq f.$  Then
  $$\lambda_0e_{\lambda_0}\leq |x| e_{\lambda_0}\leq fe_{\lambda_0},$$
i.e.
 $$\lambda_0e_{\lambda_0}\leq  fe_{\lambda_0}.$$
Thus
 $$\lambda_0 z(e_{\lambda_0})\leq  fz(e_{\lambda_0}),$$
where  $z(e_{\lambda_0})$ is the central support of the projection
$e_{\lambda_0}.$ Thus
$$\lambda_0 z(e_{\lambda_0})\leq  ||x||z(e_{\lambda_0}).$$
This means that  $||x||\neq 0.$

  Take  $g\in S(Z(M)),$ $x\in mix(M).$ We have
 $$||gx||=\inf\{f\in S(Z(M)): |gx|\leq f\}=\inf\{|g|f\in S(Z(M)): |x|\leq f\}=
 $$
 $$=|g|\inf\{f\in S(Z(M)): |x|\leq f\}=|g|||x||,$$
 i.e. $$||gx||=|g|||x||.$$

 Now let  $x, y\in mix(M).$ By  \cite[Theorem  2.4.5]{Mur} there exist partial isometries
 $u, v\in M$
such that
 $$|x+y|\leq u|x|u^{\ast}+v|y|v^{\ast}.$$
Thus
 $$|x+y|\leq u|x|u^{\ast}+v|y|v^{\ast}\leq u||x||u^{\ast}+v||y||v^{\ast}=
 $$
 $$=
 ||x||uu^{\ast}+||y||vv^{\ast}\leq||x||+||y||,$$
and therefore
$||x+y||\leq ||x||+||y||.$

Take $x, y\in mix(M).$ We may assume that
 $||x||=||y||=\textbf{1}.$
Then
 $x, y\in M$ и $||x||_{M}=||y||_{M}=1,$ and therefore
 $||xy||_{M}\leq1.$ Hence
  $|xy|\leq \textbf{1}.$
  Thus  $||xy||\leq \textbf{1},$ i.e.
   $||xy||\leq ||x||||y||.$

Let $x\in M,$   $||x||=\textbf{1}.$ By Proposition  1.3 we obtain
$$||zx||_{M}=1$$
for every  nonzero central projection
$z\in M.$
Thus
$$||zxx^{\ast}||_{M}=||(zx)(zx)^{\ast}||_{M}=
||zx||^{2}_{M}=1.$$ Therefore by Proposition
  1.3 we obtain that
$||xx^{\ast}||=\textbf{1},$
 i.e.
  $||xx^{*}||=||x||^{2}.$

Finally we shall  prove the completeness of the space $mix(M).$
First we consider the case where the center $S(Z(M))\cong
L^0(\Omega, \Sigma, \mu)$ satisfies the condition
$\mu(\Omega)<\infty.$

Let  $\{x_n\}$ be a $(bo)$-fundamental sequence in
 $mix(M),$
i.e.  $||x_n-x_m||\rightarrow0$ при $n, m\rightarrow\infty.$
 By the inequality
 $$\left| ||x_n||-||x_m|| \right|\leq ||x_n-x_m||$$
 we obtain that the sequence  $\{||x_n||\}$ is  $(o)$-fundamental in
  $S(Z(M)),$ in particular,  $\{||x_n||\}$ is order bounded in
  $S(Z(M)),$ i.e. there exists
$c\in S(Z(M))$
such that
 $||x_n||\leq c$ for all
  $n\in\mathbb{N}.$

 Now replacing $x_n$ with
 $(\textbf{1}+c)^{-1}x_n$ we may assume that
  $x_n\in M, ||x_n||\leq\textbf{1}$ and  $\{x_n\}$ is
$(bo)$-fundamental.

Since  $\mu(\Omega)<\infty$ by Egorov's theorem  for any  $k\in\mathbb{N}$ there exists
$A_k\in \Sigma$ with  $\mu(\Omega\setminus A_k)\leq\frac{\textstyle 1}{\textstyle k}$ such that
 $||\chi_{A_{k}} (x_n-x_m)||_M\rightarrow 0$ as $n, m\rightarrow\infty.$
Since  $M$ is complete, one has that $\chi_{A_{k}}x_n\rightarrow a_k$ as $n\rightarrow\infty$
 for an appropriate  $a_k\in M.$

Put
 $$
 z_1=\chi_{A_{1}},\,\,
 z_k=\chi_{A_{k}}\bigwedge(\bigvee\limits_{i=1}^{k-1}z_i)^{\perp},\,\,k\geq 2.
 $$
 Then
  $z_i\wedge z_j=0,\,\,i\neq j,$ $\bigvee\limits_{k\geq1}z_k=\textbf{1}.$
 Set  $$a=\sum\limits_{k=1}^{\infty}z_ka_k.$$ Then
 $x_n\rightarrow a.$  This means that the space  $mix(M)$ is $(bo)$-complete.

 Now we consider the  general case for the center
  $S(Z(M))\cong L^0(\Omega, \Sigma, \mu).$  There exists a mutually orthogonal system
    $\{\Omega_i: i\in I\}$ in  $\Sigma$
such that  $\mu(\Omega_i)<\infty.$ As above we have that for every
 $i\in I$
there exists  $a_i\in mix(M)$ such that
 $\chi_{\Omega_{i}}x_n\rightarrow a_i.$
 Set
  $$a=\sum\limits_{i\in I}\chi_{\Omega_{i}}a_i.$$
  Then $x_n\rightarrow a.$  This means that the space  $mix(M)$ is $(bo)$-complete.
 The proof is complete. $\blacksquare$

From Propositions  1.1 and 1.4 we obtain the following result.

\textbf{Corollary  1.5.} \emph{Let  $M$ be a  von Neumann algebra
without direct summands of type} II. \emph{Then }$(LS(M), ||\cdot||)$ \emph{is a
$C^\ast$-algebra over $S(Z(M))\cong L^{0}(\Omega, \Sigma, \mu).$}

\medskip
\begin{center} \textbf{2.   The topology  of convergence locally in measure}\end{center}

Let  $M$ be an arbitrary commutative von Neumann algebra. Then as
we have mentioned above $M$ is  *-isomorphic to the
 *-algebra $L^{\infty}(\Omega,\Sigma, \mu),$ while the algebra  $LS(M)=S(M)$ is  *-isomorphic to the
  *-algebra  $L^0(\Omega,\Sigma, \mu).$

The basis of neighborhoods of zero in the topology of convergence locally in measure
    on $L^0(\Omega,\Sigma, \mu)$ consists of the sets
$$W(A,\varepsilon,\delta)=\{f\in L^0(\Omega,\Sigma, \mu):\exists B\in \Sigma, \, B\subseteq A, \,
\mu(A\setminus B)\leq \delta, $$
$$ f\cdot \chi_B \in L^{\infty}(\Omega,\Sigma, \mu),\,
||f\cdot \chi_B||_{L^{\infty}(\Omega,\Sigma, \mu)}\leq \varepsilon\},$$
where $\varepsilon, \delta>0, \, A\in \Sigma, \, \mu(A)<+\infty.$

Recall the definition of the dimension functions on the lattice $P(M)$ of projection from $M$ (see \cite{Mur}).

By  $L_+$ we denote the set of all measurable  functions  $f: (\Omega,\Sigma, \mu)\rightarrow
[0,{\infty}]$ (modulo functions equal to zero almost everywhere ).

 Let  $M$ be an arbitrary von Neumann algebra with the center $Z=L^\infty(\Omega,\Sigma, \mu).$
Then there exists a map  $D:P(M)\rightarrow
L_{+}$ with the following properties:

(i) $d(e)$ is a finite function if only if the projection  $e$ is finite;

(ii) $d(e+q)=d(e)+d(q)$  for $p, q \in P(M),$ $eq=0;$

(iii) $d(uu^*)=d(u^*u)$ for every  partial isometry  $u\in M;$

(iv) $d(ze)=zd(e)$ for all $z\in P(Z(M)), \,\, e\in P(M);$

(v) if  $\{e_{\alpha}\}_{\alpha \in J}, \,\,\, e\in P(M) $ and $e_{\alpha}\uparrow e,$ then
$$d(e)=\sup \limits_{\alpha \in J}d(e_{\alpha}).$$
This map  $d:P(M)\rightarrow L_+,$ is a called \emph{the dimension functions} on  $P(M).$

The basis of neighborhoods of zero in the topology of convergence
locally in measure on $LS(M)$ consists (in the above notations)
of the following sets
$$V(A,\varepsilon,\delta)=\{x\in LS(M):\exists p\in P(M), \, \exists z\in P(Z(M)),  \,
xp \in M, $$
$$ ||xp||_{M}\leq \varepsilon, \,\, z^{\bot}\in W(A,\varepsilon,\delta), \,\, d(zp^{\bot})\leq \varepsilon z\},$$
where
 $\varepsilon, \delta>0, \, A\in \Sigma, \, \mu(A)<+\infty.$

We need following assertion from  \cite[pp. 242, 261, 265]{Mur}).

\textbf{Proposition  2.1.} \emph{Let
$\varepsilon, \delta>0, \, A\in \Sigma, \, \mu(A)<+\infty.$ Then:}

a) $\lambda V(A,\varepsilon,\delta)=V(A,|\lambda|\varepsilon,\delta)$
 \emph{for all} $\lambda \in \mathbb{C},\, \lambda \neq 0;$

b) $x\in V(A,\varepsilon,\delta) \Leftrightarrow |x|\in V(A,\varepsilon,\delta);$

c) $x\in V(A,\varepsilon,\delta) \Rightarrow x^{*}\in V(A,2\varepsilon,\delta);$

d) $x\in V(A,\varepsilon,\delta), \, y\in M \Rightarrow  yx\in ||y||_{M}V(A,\varepsilon,\delta);$

e) \emph{for each} $ x\in LS(M)$ \emph{there exist} $\varepsilon_1, \delta_1>0, \, B\in \Sigma, \, \mu(B)<+\infty,$
\emph{such that}
$$x\cdot V(B,\varepsilon_1,\delta_1)\subseteq V(A,\varepsilon,\delta).$$

In the next section we shall also use the following properties of the topology of convergence locally in measure.

\textbf{Lemma  2.2.} \emph{Let
 $\varepsilon, \delta>0, \, A\in \Sigma, \, \mu(A)<+\infty, \lambda \in \mathbb{C}.$
 If  $|\lambda|\leq \varepsilon,$ then} $\lambda \textbf{1}\in V(A,\varepsilon,\delta).$

Proof. Put
 $p=\textbf{1}, \, z=\textbf{1}.$ Then
  $||\lambda \textbf{1} p||_{M}=|\lambda|\leq \varepsilon,
\, z^{\bot}=0\in W(A,\varepsilon,\delta), \, D(zp^{\bot})=D(0)=0\leq \varepsilon z,$ and therefore
$\lambda \textbf{1}\in V(A,\varepsilon,\delta).$
The proof is complete. $\blacksquare$

\textbf{Lemma  2.3.} \emph{Let
 $x\in LS(M), \, \varepsilon, \delta>0, \, A\in \Sigma, \, \mu(A)<+\infty.$ Then there exists
$\lambda _{0}>0$ such that $x\in \lambda _{0}V(A,\varepsilon,\delta).$}

Proof. By Proposition  2.1 e) there exist
$\varepsilon_1, \delta_1>0, \, B\in \Sigma, \, \mu(B)<+\infty,$ such that
$$x\cdot V(B,\varepsilon_1,\delta_1)\subseteq V(A,\varepsilon,\delta).$$ From Lemma  2.2 it follows that
$\varepsilon_1 \textbf{1}\in V(B,\varepsilon_1,\delta_1).$ Therefore
$x\varepsilon_1 \textbf{1}\in V(A,\varepsilon,\delta),$ i.e.
 $x\in \lambda _{0}V(A,\varepsilon,\delta),$ where
$\lambda_0=\varepsilon^{-1}_{1}.$
The proof is complete.$\blacksquare$

\textbf{Lemma  2.4.} \emph{If  $x\in V(A,\varepsilon,\delta)$ and $u, v \in M$ are partial isometries,
then $uxv\in V(A,4\varepsilon,\delta).$}

Proof. The case when $u=0$ or  $v=0$ is  trivial. Assume that  $u, v\neq 0.$ Then
 $||u||_{M}=||v||_{M}=1.$ By Proposition  2.1 d) we obtain that
$vx\in V(A,\varepsilon,\delta).$ From Proposition  2.1 c) it follows that
 $x^*v^*=(vx)^*\in V(A,2\varepsilon,\delta).$
 Applying Proposition 2.1 d)  once more we have that
  $v^*x^*u^*\in V(A,2\varepsilon,\delta)$
and  $uxv=(v^*x^*u^*)^*\in V(A,4\varepsilon,\delta).$
The proof is complete.
$\blacksquare$

\textbf{Lemma  2.5.} \emph{If $f_i\in S(Z(M)), \, i=1,2, \, |f_1|\leq |f_2|$ and
$f_2\in V(A,\varepsilon,\delta),$ then  $f_1\in V(A,\varepsilon,\delta).$}

Proof. Let  $f_2\in V(A,\varepsilon,\delta).$ Then
 $|f_2|\in V(A,\varepsilon,\delta).$
Therefore there exist $p_0\in P(M), \, z_0\in P(Z(M))$ such that
$$|f_2|p_0\in M, \, \||f_2|p_0\|_{M} \leq \varepsilon, \, z^{\bot}_0\in W(A,\varepsilon,\delta), \,
d(z_0p^{\bot}_0)\leq \varepsilon z_0.$$
From $|f_1|\leq |f_2|$ we get
$p_0|f_1|p_0\leq p_0|f_2|p_0$ and
 $|f_1|p_0\leq |f_2|p_0.$ Hence
 $|f_1|p_0\in M$
 and
  $\||f_1|p_0\|_{M}\leq \||f_2|p_0\|_{M} \leq \varepsilon,$ i.e.
  $\||f_1|p_0\|_{M}\leq \varepsilon.$ Since
   $z^{\bot}_0\in W(A,\varepsilon,\delta), \,
d(z_0p^{\bot}_0)\leq \varepsilon z_0$ we see  that  $|f_1|\in V(A,\varepsilon,\delta)$ or
$f_1 \in V(A,\varepsilon,\delta).$
The proof is complete.
 $\blacksquare$

Recall  \cite{Sak1} that a von Neumann algebra $M$ is  said to be \emph{properly infinite},
if any nonzero central projection  $z$
in $M$ is infinite.

\textbf{Lemma 2.6.} \emph{Let  $M$ be a properly infinite von Neumann algebra and let}
$\varepsilon, \delta>0, \, A\in \Sigma, \, \delta<\mu(A)<+\infty,
\, 0<\varepsilon <1.$ \emph{Then from} $\lambda \textbf{1}\in V(A,\varepsilon,\delta),$
 \emph{where} $\lambda \in \mathbb{C},$
\emph{it follows that}  $|\lambda|\leq \varepsilon.$

Proof. Let  $\lambda \textbf{1}\in V(A,\varepsilon,\delta).$
Then there exist $p\in P(M), \, z\in P(Z(M))$ such that
 $zp^{\bot}$ is a finite and
 $\|\lambda p\|_{M}\leq \varepsilon,$
  $z^{\bot}\in W(A,\varepsilon,\delta).$
Set  $z=\chi_E,$ where
 $E\in \Sigma.$ Since $z^{\bot}\in W(A,\varepsilon,\delta)$ there exists
 $B\in \Sigma, \, B\subset A$
such that
$$\mu(A\setminus B)\leq \delta,\, \, \|z^{\bot}\chi_B\|_{M}\leq \varepsilon. $$
Since $0<\varepsilon <1,$ from the inequality  $\|z^{\perp}\chi_B\|_{M}\leq \varepsilon$ we have that
 $(\textbf{1}-\chi_E)\chi_B=0.$
From $\mu(A)>\delta$ and  $\mu(A\setminus B)\leq \delta,$ we get
 $\chi_B\neq 0,$ and therefore  from $(\textbf{1}-\chi_E)\chi_B=0$
we obtain  that
 $\chi_E\neq 0,$ i.e. $z\neq 0.$ Since  $zp^{\bot}$ is finite  and $M$ is properly infinite,
  it follows that projection $zp$ is an infinite. Therefore
 $p\neq 0.$ Thus
  $|\lambda|=|\lambda|\|p\|_{M}=\|\lambda p\|_{M}\leq \varepsilon,$ i.e.
   $|\lambda|\leq \varepsilon.$
   The proof is complete.
$\blacksquare$

\medskip

\medskip

\begin{center} \textbf{3. Additive derivations on the central extensions
of properly \\infinite von Neumann algebras}\end{center}

The following theorem is the main result of this paper.

 \textbf{Theorem 3.1.} \emph{Let  $M$ be a properly infinite  von
Neumann algebra. Then
every additive derivation on the algebra $mix(M)$ is an inner derivation.}

To prove  this theorem  we need several preliminary assertions.

Let  $\mathcal{A}$ be an algebra and denote by $Z(\mathcal{A})$ its center. If $D$ is an additive derivation
on $\mathcal{A}$ and $\Delta=D|_{Z(\mathcal{A})}$ is its restriction onto the center of $\mathcal{A},$ then
 $\Delta$ maps $Z(\mathcal{A})$ into itself \cite[Remark  1]{Alb2} (see also \cite[Lemma  4.2]{Ber1}).

Let  $M$  be a commutative  von Neumann algebra and let
$\mathcal{A}$ be an arbitrary subalgebra of $LS(M)=S(M)$ containing
$M.$ Further we shall identify the algebra $LS(M)=S(M)$ with an appropriate $L^0(\Omega,\Sigma, \mu)$.

Consider an additive derivation  $D:\mathcal{A}\rightarrow S(M)$
and let us show that  $D$ can be extended to an additive derivation
$\tilde{D}$ on the whole $S(M).$

 In the commutative von Neumann algebra $M$  for an arbitrary element  $x\in S(M)$ there exists a sequence  $\{z_n\}$
of mutually orthogonal  projections with
$\bigvee\limits_{n\in\mathbb{N}}z_n =\textbf{1}$ and  $z_n x\in M$
for all $n\in\mathbb{N}.$ Set
\begin{equation}
\tilde{D}(x)=\sum\limits_{n\geq1} z_n D(z_n x).
\end{equation}
Since every additive derivation  $D:\mathcal{A}\rightarrow S(M)$ is
identically zero on  projections of  $M,$ the equality (2)
gives a well-defined derivation
$\tilde{D}:S(M)\rightarrow S(M)$ which coincides with $D$ on
$\mathcal{A}.$

Given an arbitrary additive derivation $\Delta$ on $S(M)$=$L^0(\Omega,\Sigma, \mu)$ the element
$$z_\Delta=\inf\{\pi\in \nabla: \pi\Delta=\Delta\}$$
is called the support of the additive derivation  $\Delta,$ where $\nabla$
 is the complete Boolean algebra of all idempotents from $L^0(\Omega,\Sigma, \mu)$
 (i.e. characteristic functions of sets from $\Sigma$).

For any non trivial additive derivation $\Delta: L^0(\Omega,\Sigma, \mu)\rightarrow L^0(\Omega,\Sigma, \mu)$
 there exists a sequence
$\{\lambda_n\}_{n=1}^{\infty}$ in $L^{\infty}(\Omega, \Sigma, \mu)$ with $|\lambda_n|\leq
\textbf{1},\,n\in \mathbb{N}$, such that
$$|\Delta(\lambda_n)|\geq n z_\Delta$$ for all $n\in \mathbb{N}$ (see \cite[Lemma  2.6]{Alb2}).
In \cite{Alb2} this assertion has been proved for  linear derivations,
but the proof is the same for additive derivations.

\textbf{Lemma  3.2.} \emph{Let   $M$ be a properly infinite von
Neumann algebra, and let $\mathcal{A}\subseteq LS(M)$ be a *-subalgebra such
that  $M\subseteq \mathcal{A}$ and suppose that $D:\mathcal{A}\rightarrow \mathcal{A}$ is an
additive derivation. Then  $D|_{Z(\mathcal{A})}\equiv 0,$ in particular,
$D$ is $Z(\mathcal{A})$-linear.}

Proof.
Let  $D$ be an additive  derivation on $\mathcal{A},$ and let $\Delta$
be its restriction onto   $Z(\mathcal{A}).$ Since $M\subset \mathcal{A}\subset LS(M)$  it follows that
$Z(M)\subset Z(\mathcal{A})\subset S(Z(M))$ = $L^0(\Omega,\Sigma, \mu)$.
Let us extend the additive derivation
$\Delta$  onto whole $S(Z(M))$ as in (2) above, and denote the extension also
 by $\Delta.$

Since  $M$ is  properly infinite there exists a sequence of
mutually orthogonal  projections
$\{p_n\}_{n=1}^{\infty}$ in $M$ such that  $p_n\sim\textbf{1}$ for
all $n\in\mathbb{N}$, and
  $\bigvee\limits_{n=1}^{\infty}p_n=\textbf{1}.$

For any bounded sequence  $\Lambda=\{\lambda_n\}_{n\in\mathbb{N}}$ in $Z(M)$
define an operator  $x_\Lambda$ by
$$x_\Lambda=
\sum\limits_{n=1}^{\infty}\lambda_n p_n.$$ Then
\begin{equation}
x_\Lambda p_n=p_n x_\Lambda=\lambda_n p_n
\end{equation}
for all $n\in\mathbb{N}.$

Take  $\lambda\in Z(\mathcal{A})$ and $n\in \mathbb{N}.$ From the identity
$$D(\lambda p_n)=D(\lambda)p_n+\lambda D(p_n)$$
 multiplying it by $p_n$ from the both sides we obtain
$$p_nD(\lambda p_n)p_n=p_n D(\lambda)p_n+\lambda p_n D(p_n)p_n.$$
Since  $p_n$ is a projection, one has that  $p_n D(p_n)p_n=0,$ and
since $D(\lambda)=\Delta(\lambda)\in Z(\mathcal{A}),$ we have
\begin{equation}
p_nD(\lambda p_n)p_n=\Delta(\lambda)p_n.
\end{equation}

Now from the identity
$$D(x_\Lambda p_n)=D(x_\Lambda)p_n+x_\Lambda D(p_n),$$
in view of (3)  one has similarly
$$p_nD(\lambda_n p_n)p_n=p_n
D(x_\Lambda)p_n+\lambda_n p_n D(p_n)p_n,$$ i.e.
\begin{equation}
p_n D(\lambda_n p_n)p_n=p_nD(x_\Lambda)p_n.
\end{equation}
Now (4) and (5) imply
 \begin{equation}
p_nD(x_\Lambda)p_n=\Delta(\lambda_n)p_n.
\end{equation}

If we suppose that $\Delta\neq 0$ then  $z_\Delta\neq0.$ By \cite[Lemma 2.6]{Alb2} there exists a
bounded sequence $\Lambda=\{\lambda_n\}_{n\in\mathbb{N}}$ in $Z(M)$ such that

\medskip

$$|\Delta(\lambda_n)|\geq n  z_\Delta$$
for all $n\in \mathbb{N}.$

Replacing the algebra  $M$ by the algebra  $z_\Delta M,$ and  the additive derivation
 $D$ by $z_\Delta D$,  we may assume that
$z_\Delta=\textbf{1},$ i.e.
\begin{equation}
|\Delta(\lambda_n)|\geq n \textbf{1}
\end{equation}
 for all  $n\in \mathbb{N}.$

Now take  $\varepsilon, \delta>0,$
$A\in\Sigma,$ $\delta< \mu(A)<+\infty.$ By Lemma  2.3 there exists a number
$\lambda_0>0$  such that  $D(x_\Lambda)\in\lambda_0 V(A, \varepsilon, \delta).$
From Lemma  2.4 it follows that    $p_nD(x_\Lambda)p_n\in\lambda_0 V(A, 4\varepsilon, \delta)$
 for all $n\in \mathbb{N}.$
If we combine this with   (6) we obtain
\begin{equation}
\Delta(\lambda_n)p_n\in\lambda_0 V(A, 4\varepsilon, \delta)
\end{equation}
for all
 $n\in \mathbb{N}.$  Since  $p_n\sim\textbf{1}$  for each
 $n\in\mathbb{N}$ , there exists a sequence of partial isometries  $\{u_n\}_{n\in\mathbb{N}}$ in $M$ such that
$u_nu_n^{*}=p_n$  and $u_n^{*}u_n=\textbf{1}$ for all
$n\in\mathbb{N}.$ Using (8) and Lemma  2.4 we have
$$
u_n^{*}\Delta(\lambda_n)p_nu_n\in\lambda_0 V(A, 16\varepsilon, \delta)$$
for all $n\in \mathbb{N}.$ Thus from the  equality
$u_n^{*}p_nu_n=u_n^{*}u_nu_n^{*}u_n=\textbf{1},$
we obtain that
$$
\Delta(\lambda_n)\in\lambda_0 V(A, 16\varepsilon, \delta)$$
for all $n\in \mathbb{N}.$
 Thus by Lemma  2.5 and from the  inequality  (7) we have
 \begin{equation}
 n\textbf{1}\in\lambda_0 V(A, 16\varepsilon, \delta)
\end{equation}
for all $n\in \mathbb{N}.$
Take a number  $n_0\in\mathbb{N}$ such that  $n_0>16\lambda_0\varepsilon.$ From
Proposition 2.1 a) and  (9) we obtain that
 $$
 \textbf{1}\in V(A, 16\lambda_0 \varepsilon n_{0}^{-1}, \delta).
$$
Since  $\delta< \mu(A)$ and $16\lambda_0 \varepsilon n_{0}^{-1}<1,$
from Lemma  2.6 we have the inequality  $1\leq 16\lambda_0 \varepsilon n_{0}^{-1},$ which contradicts
the inequality
$n_0>16\lambda_0\varepsilon.$ This contradiction implies that
  $\Delta\equiv 0,$ i.e. $D$ is identically
zero on the center of $\mathcal{A},$ and therefore it is $Z(\mathcal{A})$-linear.
The proof is complete. $\blacksquare$

\textbf{Remark.} A result similar to Lemma 3.2 for the case of linear derivations has been announced without proof in \cite[Proposition 6.22]{Ber1}.

In the case of linear derivations on the algebras  $\mathcal{A}=S(M)$ or $S(M, \tau)$
a shorter proof of Lemma   3.2 can be obtained also from the following
result.

 \textbf{Proposition 3.3.} \emph{Let  $M$ be a properly infinite  von
Neumann algebra with the center $Z(M).$ Then
 the centers of the algebras $S(M)$ and  $S(M,\tau)$ coincide with $Z(M).$}

   Proof. Consider a central element  $z\in S(M), \ z\geq 0,$   and let
     $z=\int\limits_{0}^{\infty}\lambda\,de_{\lambda}$ be its spectral resolution.
      Then  $e_{\lambda}\in Z(M)$ for all
     $\lambda>0.$ Assume that  $e_{n}^{\bot}\neq 0$ for all  $n\in\mathbb{N}.$
     Since $M$ is properly infinite, $Z(M)$ does not contain non-zero finite projections. Thus
   $e_{n}^{\bot}$ is infinite for all  $n\in\mathbb{N},$ which
   contradicts  the condition  $z\in S(M).$ Therefore there
   exists  $n_{0}\in\mathbb{N}$ such that  $e_{n}^{\bot}=0$ for all  $n\geq n_0,$
  i.e. $z\leq n_0\textbf{1}.$ This means that   $z\in Z(M),$ i.e. $Z(S(M))=Z(M).$ Similarly $Z(S(M, \tau))=Z(M).$
  The proof is complete. $\blacksquare$

  Let  $M$ be a properly infinite  von
Neumann algebra with the center $Z(M)$ and let
  $D$ be a linear   derivation on the algebra  $\mathcal{A}=S(M)$ or $S(M, \tau).$
    Proposition  3.3 implies  that  $Z(\mathcal{A})=Z(M),$
and therefore  $\Delta=D|_{Z(\mathcal{A})}$ is a linear derivation on the algebra  $Z(M).$
By  \cite[Lemma  4.1.2]{Sak1} we obtain that  $\Delta=0$ as it was asserted in Lemma 3.2.

\textbf{Proof  of Theorem 3.1.}  Let  $D: mix(M)\rightarrow mix(M)$ be an additive derivation.
From Lemma  3.2 it follows that $D$ is $S(Z(M))$-linear.
By Proposition  1.4  we have that $mix(M)$ is a $C^{\ast}$-algebra over  $S(Z((M))\cong L^0(\Omega,
\Sigma, \mu).$
Since  $D$ is  $S(Z(M))$-linear, by  \cite[Theorem  5]{AK} we obtain that  $D$ is a $S(Z(M))$-bounded,
 i.e. there exists $c\in S(Z(M))$ such that  $||D(x)||\leq c||x||$
for all $x\in mix(M).$ Take
 a sequence of
pairwise  orthogonal  central projections
$\{z_n\}_{n\in\mathbb{N}}$ in  $M$ with $\bigvee\limits_{n\geq1}z_n=\textbf{1}$ such that
$z_n c\in Z(M)$ for all $n.$ Then for any
$x\in M$ we have
$$||D(z_n x)||=z_n||D(x)||\leq z_n c ||x||,$$
i.e. $||D(z_n x)||\in Z(M).$ Thus
$$z_nD(x)\in z_n M. $$ Therefore the operator   $z_n D$ maps each subalgebra  $z_n M$ into itself for all
 $n\in\mathbb{N}.$
By Sakai's  theorem \cite[Theorem  4.1.6]{Sak1}  there exists  $a_n\in z_n M$ such that
$$
z_nD(x)=a_nx-xa_n,\,\,x\in z_n M.
$$ Set  $a=\sum\limits_{n\geq1}z_n a_n.$ Then
$a\in mix(M)$ and $D(x)=ax-xa$ for all $x\in mix(M).$ This means that
 $D$ is inner. The proof is complete. $\blacksquare$

From Theorem  3.1 and Proposition  1.1 we obtain the following result
which generalizes and extends  Theorem
 2.7 from \cite{Alb2}.

\textbf{Corollary  3.4.} \emph{Let  $M$ be a  direct sum of von Neumann algebras
of type} I$_\infty$ \emph{and} III. \emph{Then every additive derivation on
the algebra $LS(M)$ is an inner derivation.}

 \newpage

\textbf{Acknowledgments.} \emph{Part of this work was done within
the framework of the Associateship Scheme of the Abdus Salam
International Centre  for Theoretical Physics (ICTP), Trieste,
Italy. The first author thank ICTP for providing financial support
 and all facilities (July-August, 2009).}

 Sh.A.Ayupov (corresponding author):

\emph{Institute of Mathematics and Information  Technologies,}

\emph{Uzbekistan Academy of Sciences, Tashkent, Uzbekistan}

and

\emph{Abdus Salam International Centre for Theoretical Physics, Trieste, Italy,}

 permanent address:  \emph{Dormon Yoli str. 29, 100125, Tashkent, Uzbekistan}

 e-mail: \emph{sh\_ayupov@mail.ru}

\medskip
K.K. Kudaybergenov:

\emph{ Karakalpak state university,}

address: \emph{ Ch. Abdirov str. 1, 142012, Nukus, Uzbekistan,}

 e-mail: \emph{karim2006@mail.ru}

\end{document}